\documentclass{amsart}

\usepackage{amsthm}
\usepackage{graphicx}
\usepackage{amssymb}
\newtheorem{teo}{Theorem}[section]

\newtheorem{prop}[teo]{Proposition}
\newtheorem{cor}[teo]{Corollary}

\newtheorem{question}[teo]{Question}

\theoremstyle{definition}
\newtheorem{defi}[teo]{Definition}
\newtheorem{rem}[teo]{Remark}
\newtheorem{example}[teo]{Example}

\theoremstyle{remark}
\newtheorem{prof}[teo]{Proof of}

\def\mr{\mathbb{R}}

\begin{document}
\title{A short introduction to shadows of $4$-manifolds}

\author[Costantino]{Francesco Costantino}
\address{Institut de Recherche Math\'ematique Avanc\'ee\\
  Rue Ren\'e Descartes 7\\
  67084 Strasbourg, France}
\email{f.costantino@sns.it}


\begin{abstract}
We give a self-contained introduction to the theory of shadows as
a tool to study smooth $3$ and $4$-manifolds. The goal of the
present paper twofold: on one side it is intended to be a shortcut
to a basic use of the theory of shadows, on the other side it
gives a sketchy overview of some of the most recent results on shadows. No original result is proved here
and most of the details of the proofs are not included.
\end{abstract}

\maketitle

\tableofcontents
\section{Introduction}
Shadows were defined by V. Turaev for the first time at the
beginning of the nineties in \cite{Tu3} as a method for
representing knots alternative to the standard one based on knot
diagrams and Reidemeister moves. The theory was then developed in
the preprint ``Topology of shadows" which was later included in a
revised form in \cite{Tu}; moreover, a short account of the
theory was published by Turaev in \cite{Tu2}. Since then, probably
due to the slightly higher degree of complication of this theory
with respect to Kirby calculus, only few applications of shadows
were studied. Among these applications we recall the use of
shadows to study Jones-Vassiliev invariants of knots made by U.
Burri in \cite{Bu} and A. Shumakovitch in \cite{Sh} and the study
of ``Interdependent modifications of links and invariants of
finite degree" developed by M.N. Goussarov in \cite{Gou}.

In the last two years, there has been a new surge of interest in
shadows testified by the appearance of the preprint \cite{Th} of
D. Thurston on knotted trivalent graphs and shadows and by the
study of a new notion of complexity of $3$ manifolds based on
shadows developed by the author and D. Thurston in \cite{CT2}.
This notion of ``shadow complexity" turns out to be intimately
connected with the geometric properties of $3$-manifolds and in
particular with their hyperbolic structures; this fact was also
used in \cite{CFMP} to exhibit and study a particular class of
$3$-manifolds having strong geometric properties. Moreover a
refinement of the notion of shadow of $4$-manifold, called
``branched shadow" has been studied by the author in \cite{Co4},
to attack a combinatorial study of gauge invariants and complex
convexity problems. Branched shadows allow one to encode homotopy
classes of almost complex structures on $4$-manifolds, which, on
manifolds admitting a shadow, are in a natural bijection with the
set of $Spin^c$-structures. Moreover, it has been proved by the
author that branched shadows embedded in complex manifolds behave
like real surfaces do; in particular a classification of the
points of a branched shadow into totally real, complex elliptic
and complex hyperbolic ones is possible. These arguments were then
used to study geometric problems as the existence of
integrable representatives of homotopic classes of almost complex
structures or the existence of Stein domain structures on the
neighborhood of a shadow.

Roughly speaking, shadows of $4$-manifolds are $2$-dimensional
polyhedra embedded in $4$-manifolds equipped with extra-structures
suitable to encode the topology of their regular neighborhood;
they represent the analogous of spines of $3$-manifolds in the
$4$-dimensional case. Shadows allow a completely combinatorial
approach similar to the one used in the $3$-dimensional
case by means of spines.

Since the first obstacle to the comprehension of shadows is the
familiarity with simple polyhedra, we inserted a first section in
which we clarify what a simple polyhedron is and how it can be
combinatorially constructed by means of a finite set of basic
building blocks. Sections \ref{sec2},\ref{sec3} and \ref{sec4} are
devoted to expose the basic facts of the theory of shadows proved
by Turaev in \cite{Tu}. In Section \ref{sec2}, we study how an
embedding of a polyhedron into a $4$-manifold can be used to equip
the polyhedron with an additional structure, called ``gleam"
by Turaev, which is fundamental in the theory. In Section
\ref{sec3}, we recall Turaev's method to thicken a polyhedron
equipped with gleams and to construct $4$-manifolds. In Section
\ref{sec4}, we introduce a set of basic moves on shadows and study
their effect on the corresponding $4$-manifolds.

The last three sections are devoted to an overview of some of
the open problems and questions. 
In Section \ref{sec5}, we outline connections between shadows of
$4$-manifolds and the Andrews-Curtis conjecture. Section
\ref{sec6} is devoted to the results proven
recently by the author and D. Thurston on shadow complexity and
hyperbolic geometry of $3$-manifolds. At last, in Section
\ref{sec7} we outline some of the recent developments of
the theory towards a combinatorial method for studying
$4$-manifolds equipped with additional structures such as
$Spin^c$, almost-complex, complex and Stein structures.

{\bf Acknowledgements.} The author wishes to warmly thank Riccardo
Benedetti, Dylan Thurston and Vladimir Turaev for their criticism
and encouraging comments.

\section{Simple polyhedra}\label{sec1}
\begin{figure}
  \centerline{\includegraphics[width=8.4cm]{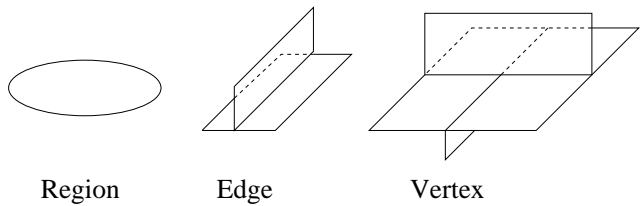}}
  \caption{The local models of a simple polyhedron.  }
  \label{fig:singularityinspine}
\end{figure}

The basic objects underlying shadows of $4$-manifolds are simple
$2$-dimensional polyhedra, that is $2$-dimensional CW-complexes
whose points have regular neighborhoods of ``simple" type as those
shown in Figure \ref{fig:singularityinspine}.
From the figure, we observe that the set of points whose regular neighborhood 
in a simple polyhedron $P$ is shaped on the
last two models, and which is called the {\it Singular set} of the
polyhedron, is a $4$-valent graph denoted by $Sing(P)$ whose
vertices are exactly the points whose regular neighborhoods are
shaped on the third model and are called {\it vertices} of the
polyhedron. The complement of the singular set of the polyhedron
is the union of a set of open surfaces embedded in the
polyhedron and called {\it regions}.

One check that the local
pattern around a vertex is homeomorphic to the regular
neighborhood of the vertex of the cone over the edges of a
tetrahedron; this shows that the full symmetry group
of the regular neighborhood of a vertex is the group of
permutations of $4$ elements. Another way of
visualizing a vertex is to consider the result of attaching the
polyhedron $Y\times [0,1]$ (where $Y$ is a $Y$-shaped graph) along
$Y\times {0}$ on a $Y$-shaped graph embedded in a disc.

To clarify the reason why we restrict only to the three kind of
local singularities of Figure \ref{fig:singularityinspine} we
shift to the one dimensional case where the analogous of simple
polyhedra is the set of graphs whose local possible patterns are the one 
homeomorphic to the open interval and the one
homeomorphic to an open $Y$-shaped graph; in this case the
vertices are the vertices of the graph. A generic graph could
contain some $n$ (with $n>3$) valent vertex, so that the
restriction to $3$-valent graphs appears somewhat strong. Roughly
speaking, the idea underlying this choice is that a $n$-valent
vertex of a graph is ``stable" only if $n=3$: if we
slightly perturb the position of an edge of the graph near, say, a
$4$-valent vertex, we can transform by a local modification
the vertex into a pair of $3$-valent vertices, passing from a
$X$-shaped graph to a $H$-shaped one; if we try to repeat this
operation for one of the two $3$-valent vertices obtained that
way, we cannot further modify the homeomorphism type of the graph.
More formally (but we will not plug into these matters), each
graph is simply-homotopic to a trivalent one. Analogously, in
dimension $2$, each polyhedron can be suitably perturbed to
a new one whose singularities are those depicted in Figure
\ref{fig:singularityinspine}.

The first example of simple polyhedron is a closed
surface, whose points are all of the first type in Figure
\ref{fig:singularityinspine}. An example containing
only points of the first two types is the union of a sphere in $\mr^3$ and the disc spanned by a maximal
circle in it and contained in the ball bounded by the sphere. This
polyhedron has singular set composed by a circle and contains no
vertices; its regions are three (open) discs.

A third example can be constructed by using the last one and
gluing to it along another equator (different from the one already
chosen) another disc, this time contained in the complement of the
ball bounded by the initial sphere. The polyhedron we get has a
singular set containing two $4$-valent vertices 
(the intersections of the two maximal circles in the sphere
used to construct it) and such that no edge of the graph has
endpoints in the same vertex; the regions of the polyhedron are
six discs.

In the above example we used a sphere embedded in
$\mr^3$ but this was only for the sake of clarity: one the main
properties of simple polyhedra is that they can be abstractly
described by means of their combinatoric structure. To give
an idea of how this can be performed we do this for the preceding
example.

First, one identifies the structure of the singular set: in
our case a $4$-valent graph with two vertices connected by
four edges. Then, one studies how each region is glued to the
singular graph: in our case the boundary of each region
passes exactly once over two different edges of the Singular set
and zero times over the others. Moreover, for each pair of edges
there is exactly one region whose boundary is the union of the two
edges: there are $4$ edges and hence $6$ pairs of edges and
exactly $6$ regions. Hence, to reconstruct the polyhedron, 
it is sufficient to start from the $4$-valent
graph and glue the six regions over all the possible pairs of
edges.

In the preceding case, we were treating a particular type of
polyhedron: all its regions were discs and its singular graph was
connected; such polyhedra are called {\it standard}. The
above combinatorial reconstruction is applicable also in
the general case of simple-polyhedra, but it is easier to
describe in the standard case since there is no need of specifying 
the topology of each region to be glued to the singular graph. From
now on, for the sake of brevity, we will use the word polyhedron
for standard polyhedron, restricting the set of objects we
will study; anyway, all our discussions will hold also for simple
polyhedra in general.

Before ending this section, we define the notion of ``mod $2$
gleam" of a polyhedron $P$. Let $D$ be a region of $P$, and let us
consider the regular neighborhood of $\partial D$ in $P$. For the
sake of simplicity, we suppose that $D$ is embedded in $P$, but the
same arguments apply to the general case when $\partial D$ runs
more than once over an edge of $Sing(P)$. The regular neighborhood
of a point $p_0$ of $\partial D$ different from a vertex of $P$
is composed by a half plane (corresponding to $D$) glued
along a segment to a square in $\mr^2$ which is split in two
components, corresponding to the other regions touching $D$ along
its boundary near $p_0$. Let us choose one of these 
components, and continue this choice by following $\partial
D$ until we reach a vertex of $P$. 
To extend this choice near the vertex, we
delete from the neighborhood of the vertex the only region which
locally touches $D$ only in the vertex itself (there are six
regions around a vertex and they are ``opposite" in pairs): we
obtain again a square divided in two by $\partial D$ and we
can continue our choice canonically (we recommend the reader to
try an visualize it). After completing a whole loop around
$\partial D$, we come back to the initial point $p_0$ with a choice
of one of the two components into which $\partial D$ splits the
neighborhood of $p_0$ in $P- D$. If this second component is 
the one we choose initially then we say that the mod $2$
gleam of $D$ is $0$, otherwise it is $1$. Acting similarly for
each region of $P$ we get a $\mathbb{Z}_2$-coloring on the regions
of $P$ called {\it mod $2$ gleam} of $P$.

In the case of an embedded region $D$, the above description 
can be summarized as follows: consider the regular neighborhood of $\partial
D$ in $P- int(D)$; if it collapses over a Moebius strip, then the
mod $2$ gleam of $D$ is $1$, otherwise it is $0$.

\section{Polyhedra in $4$-manifolds}\label{sec2}

In this section we investigate how a polyhedron $P$ embedded in a
$4$-manifold $M$ can be given an extra structure related to the
topology of its regular neighborhood. Let us be more specific regarding the word
``embedded":
\begin{defi}
A polyhedron embedded in a smooth $4$-manifold is said to be {\it
locally flat} if for each point $p$ of $P$ there is a local chart
$(U,\phi)$ of the smooth atlas of $M$ such that the image of
$P\cap U$ through $\phi$ is contained in a $3$-dimensional plane
in which, up to a self-homeomorphism of the plane, it appears
as one of the three local pictures of Figure
\ref{fig:singularityinspine} in $\mr^3\subset \mr^4$, that is,
around each of its points, $P$ is contained in a $3$-dimensional
slice of $M$ and in this slice it appears (up to homeo) as shown
in Figure \ref{fig:singularityinspine}.
\end{defi}

\begin{figure}
  \centerline{\includegraphics[width=8.4cm]{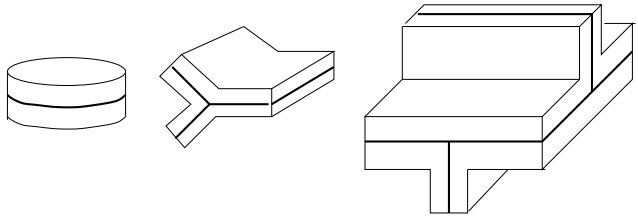}}
  \caption{The three type of blocks used to thicken a spine of a $3$-manifold.  }
  \label{fig:spineblocks}
\end{figure}

For the sake of brevity, from now on we will use the word
``embedded" for ``locally flat embedded". The first question we
ask ourselves is the following: can we reconstruct the regular
neighborhood in a manifold $M$ of a polyhedron $P$ from its
combinatorics? Before trying and ask this question in the
$4$-dimensional case, let us understand the easier $3$-dimensional
case, where $P$ is a polyhedron embedded in an oriented
$3$-manifold $N$. In this case one sees that the
$\mathbb{Z}_2$-gleam of each region of $P$ is $0$: indeed, 
roughly speaking, the regular neighborhood of a region of $P$
(which is a disc) in $N$ is $D^2\times [- 1,1]$ and
so it is disconnected by the region; this would be false if the
mod $2$ gleam of the region was $1$. Hence,
the fact that a polyhedron embeds in a $3$-manifold imposes some
restrictions on its combinatorics. The truth is much deeper: the
combinatorics of the polyhedron allows one to reconstruct the
topology of its regular neighborhood in $N$. Indeed, by
decomposing $P$ in the local patterns shown in Figure
\ref{fig:singularityinspine}, we can decompose its regular
neighborhood in blocks as those of Figure
\ref{fig:spineblocks} which can be reglued to each other
according to the combinatorics of $P$. That way, a polyhedron with
zero mod $2$ gleam determines a $3$-manifold
collapsing onto it. It can be shown that any $3$-manifold with
non-empty boundary can be constructed that way, and so the study
of three manifolds can be attacked by means of particular
polyhedra; this is the basic idea underlying the theory of spines
of $3$-manifolds.

Let us now skip to the $4$-dimensional case and, for the sake of simplicity, choose 
$P$ to be homeomorphic to $S^2$.
Suppose that $P$ is embedded in a smooth, oriented
$4$-manifold $M$: it is possible to reconstruct its regular
neighborhood using only its topology? The answer
is ``no" since the regular neighborhood of a sphere (and more
generally of a surface) is determined by the topology of the
surface and by its self-intersection number in the manifold. To
state it differently, the regular neighborhood of a surface in a
$4$-manifold is diffeomorphic to the total space of a disc bundle
over the surface (its normal bundle) and the Euler number of this
bundle is a necessary datum to reconstruct its topology.

Hence, we see that to codify the topology of the
regular neighborhood of $P$ in $M$ we need to decorate $P$ with
some additional information; when $P$ is a surface, 
the Euler number of its normal bundle is a sufficient datum.

More generally the following holds:
\begin{prop}\label{embpoly}
Let $P$ be a polyhedron embedded in an oriented smooth
$4$-manifold $M$. There is a coloring of the regions of $P$ with
values in the half integers $\frac{1}{2}\mathbb{Z}$ canonically
induced by its embedding in $M$. This coloring is called the {\it
gleam} of $P$ in $M$.
\end{prop}

\begin{prf}{1}{
The construction we will describe is similar to the one
used to define the mod $2$ gleam of a polyhedron. For the sake of
simplicity let again $D$ be a region of $P$ whose boundary is
embedded in $Sing(P)$. Let $p_0$ be a point of $\partial D$ and
let us consider the three dimensional slice $B^3$ of $M$ into
which a regular neighborhood of $p_0$ in $P$ is sitting (it exists
since by hypothesis $P$ is locally flat in $M$). Let us fix a
auxiliary riemannian metric on $M$ and consider the
orthogonal direction to $D$ in $B^3$; by construction, 
it coincides with the direction along which the two other
regions touching $\partial D$ on $p_0$ get separated (see Figure
\ref{fig:divergingdirection}).
\begin{figure}
  \centerline{\includegraphics[width=6.4cm]{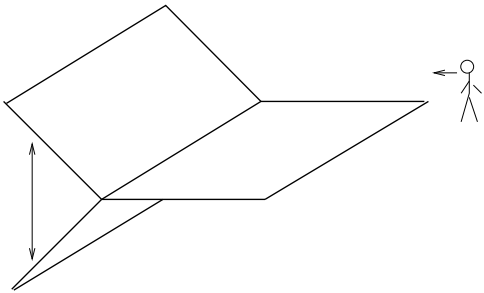}}
  \caption{The picture sketches the position of the polyhedron in a $3$-dimensional slice of the ambient $4$-manifold.
  The direction indicated by the vertical double arrow is the one
  along which the two regions touching the horizontal one get separated.  }
  \label{fig:divergingdirection}
\end{figure}

We can extend the definition of this direction to the whole
$\partial D$, in the same way we did to get the mod $2$ gleam,
obtaining a continuous choice of orthogonal directions to $D$ in
$M$ along $\partial D$. This choice represents a section of 
the bundle of orthogonal directions to $D$ in $M$, whose
fiber is $S^1$. The obstruction to extend this section to the
whole $D$ can be canonically identified with an integer: indeed it
represents an element of $H^2(D,\partial D;\pi_1(S^1))$ and if we
fix an orientation on $D$ then also the fiber gets oriented
since $M$ is oriented, hence giving an identification of
$H^2(D,\partial D;\pi_1(S^1))$ with $H^2(D,\partial
D;\mathbb{Z})=\mathbb{Z}$; this identification does not change if
we revert the orientation of $D$, since also the orientation of
the fiber changes. Hence, if we divide by $2$ the number obtained
through this identification, we get an element of
$\frac{1}{2}\mathbb{Z}$ which we will call the gleam of $D$, and
which represents the Euler number of the normal bundle of $D$ in
$M$. It is worth note that the Euler number of the normal bundle
of an embedded surface equals the self-intersection number of
the surface, so that the gleam represents a version of the
self intersection number  ``localized" on the regions of
the polyhedron.

The above construction can be performed analogously in the
case when $D$ is not embedded in $P$; that way, we can define a
gleam over each region of $P$. Note that, by construction, the
gleam of a region of $P$ is non-integer (an odd multiple of
$\frac{1}{2}$) if and only if its mod $2$ gleam is $1$.

}\end{prf}

\section{From polyhedra with gleams to $4$-manifolds}\label{sec3}

In the preceding section we showed that any polyhedron embedded in
a $4$-manifold can be equipped with a gleam. This produces a map
from the pairs $(M,P)$ of smooth, oriented and compact manifolds
containing embedded polyhedra into the set of polyhedra equipped
with gleams $(P,gl)$ where a gleam is a $\frac{1}{2}\mathbb{Z}$
coloring of the regions of $P$ such that the color of a region is
non integer if and only if its mod $2$ gleam is $1$. We also
stated that if $P$ is a surface then its gleam is 
the Euler number $e$ of the normal bundle of the surface in
$M$ and hence the pair $(P,e)$ is sufficient to reconstruct the
regular neighborhood of $P$ in $M$. This is true in general as
proved by Turaev in \cite{Tu}:
\begin{teo}[Reconstruction map]\label{teo:reconstruction}
Let $(P,gl)$ be a polyhedron equipped with gleams; there exists a
canonical reconstruction map associating to $(P,gl)$ a pair
$(M_P,P)$ where $M_P$ is a smooth, compact and oriented
$4$-manifold containing an embedded copy of $P$ over which it
collapses and such that the gleam of $P$ in $M_P$ coincides with
the initial gleam $gl$. The pair $(M_P,P)$ can be explicitly
reconstructed from the combinatorics of $P$ and from $gl$.
Moreover, if $P$ is a polyhedron embedded in a smooth and oriented
manifold $M$ and $gl$ is the gleam induced on $P$ as explained in
the preceding section, then $M_P$ is diffeomorphic to a regular
neighborhood of $P$ in $M$.
\end{teo}

The above theorem is the key tool of the theory of shadows of
$4$-manifolds: indeed by means of this result, one is allowed to
study $4$-manifolds in a purely combinatorial way, using only
polyhedra equipped with gleams.
\begin{prf}{1}{
The idea of the proof is to decompose $P$ in blocks as those of
Figure \ref{fig:singularityinspine}, thicken them to canonical
blocks, and then describe how to glue these blocks according to
the combinatorics of $P$ and to its gleam. The shape of the basic
blocks is easily understood: consider the blocks used to
reconstruct $3$-manifolds from their spines and take their
products with $[- 1,1]$. Then the block corresponding to a region
is $D^2\times D^2$ and the blocks corresponding to edges and
vertices of $Sing(P)$ are the products respectively of the second
and third block of Figure \ref{fig:spineblocks} with a segment
$[- 1,1]$ representing ``the fourth dimension".

We first reconstruct the regular neighborhood in $M_P$ of the
subpolyhedron $S$ of $P$ coinciding with the regular neighborhood
of $Sing(P)$ in $P$; to do this, we use the blocks corresponding
to vertices and edges of $P$. The combinatorics of $P$ forces us
to glue these blocks to each other and leaves us only one choice:
whether gluing them by reverting or not the $[- 1,1]$ factor of the
fourth dimension of the blocks. If $P$ contains
$n$ vertices, we have to connect the $n$ blocks
corresponding to these vertices through the $2n$ blocks
corresponding to the edges of the $4$-valent graph $Sing(P)$.
Consider a maximal subtree $T$ of $Sing(P)$, and glue the blocks
corresponding to the vertices and the edges of $T$ according to
the combinatorics of $S$ choosing arbitrarily the gluing on
the level of the $[- 1,1]$ factors, that is choosing arbitrarily
wether identifying the $[- 1,1]$ fibers by the identity or the
multiplication by $- 1$. The
manifold constructed that way is $B^4$ and contains a properly
embedded copy of a regular neighborhood of $T$ in $S$. Only $n+1$
blocks corresponding to the edges of $Sing(P)- T$ are now left to
be glued according to the combinatorics of $S$. For each of these
blocks, there is only one way of gluing the $[- 1,1]$ fiber
to get an orientable $4$-manifold; moreover, 
at the end of this construction, the manifold $H$ obtained 
is diffeomorphic to the regular neighborhood of a $4$-valent graph
($Sing(P)$) in $\mr^4$.

Note
that since $H$ admits an orientation reversing diffeomorphism
corresponding to the product by $- 1$ along the fiber $[- 1,1]$ of
the blocks we used to construct it, then it is canonically
oriented. Moreover, by construction, $H$ fibers over a (not
necessarily orientable) $3$-manifold $L$ whose spine is $S$; the
fiber is exactly the $[- 1,1]$ factor of the blocks used to
construct $H$. To construct $L$ 
it is sufficient to repeat the above construction of $H$ using the
$3$-dimensional blocks of Figure \ref{fig:spineblocks} 
without taking the product with $[- 1,1]$; $L$ is not
necessarily orientable and, since its construction does not
deal with the regions of $P$, no assumption on the gleams of the
regions of $P$ is necessary to construct $L$ or $H$. Finally, a
copy of $S$ is properly embedded in $L$; then it is
embedded in a locally flat way in $H$ so that its boundary is contained
in $\partial H$. Each boundary component of $S$
corresponds by construction to the boundary component of a region
of $P$, hence to finish our construction we have to glue
along $\partial S\subset
\partial H$ the blocks corresponding to the regions of $P$.

Let us show how to do this for a single region $Y$ of $P$. The block
corresponding to $Y$ is $Y\times D^2$ and, since by
hypothesis $Y$ is a disc, it is a $4$-ball with a prescribed
$2$-handle structure, $Y$ being the core of the handle. 
To describe how to glue this $2$-handle to $H$ we:
\begin{enumerate}
\item identify the image of $\partial Y$ in $\partial H$
through the gluing map;

\item identify the image in $\partial H$ through the gluing map
of the canonical framing of $\partial Y$, that is of the curve
$\partial Y\times \{1\}$.
\end{enumerate}

The first step is straightforward: $\partial Y$ has to be
identified with the boundary component $s$ of $S$ contained in
$\partial H$ corresponding to $Y$. To perform the second step, we
have to identify a canonical framing on $s$; this can be done
since $H$ fibers over $L$ with fiber $[- 1,1]$. Indeed, while assembling
$H$, we glued blocks having the structure of products
between $3$-dimensional blocks and $[- 1,1]$ and we always
respected this product structure up to the orientation of the
$[- 1,1]$-fiber. Hence each boundary component of $S\subset H$ (and
in particular $s$) is automatically equipped with the framing
given by the $[- 1,1]$-direction; the only delicate point is that
the fiber over $s$ of the projection of $H$ onto the $L$ could in some cases 
be a Moebius strip instead of an annulus.. If it is an
annulus, we have a genuine framing on $s$ and hence, 
to attach the $2$-handle corresponding to $Y$, is
sufficient to say how many times the image of the canonical
framing of $\partial Y$ in $Y\times D^2$ twists around $s$ in
$\partial H$ with respect to this framing; since $H$ is oriented,
an integer number is then sufficient to encode this information, 
and we use the gleam of $Y$. Indeed, the gleam of $Y$ is
integer exactly when the inverse image in $H$ of the curve $s$
through the projection along the $[- 1,1]$-fiber is an annulus.

We are left with the case when the inverse image of $s$ in $H$
through the projection along the $[- 1,1]$-direction is a Moebius
strip $M_s$. In this case, the gleam of $Y$ is half integer; we
then add a positive half twist to $M_s$ around $s$ (recall that
$H$ is oriented), obtaining an annulus whose core is $s$ and hence
a framing on $s$. We now perform the construction of the preceding
case using as new gleam for $Y$ the integer given by the initial
gleam of $Y$ minus $\frac{1}{2}$.

Performing the above construction for all the regions of $P$ we
get a $4$-dimensional ``neighborhood" of $P$; it can be checked
that the gleam induced on $P$ by this neighborhood coincides with
the initial gleam $gl$.}\end{prf}

\begin{example}
Let $P$ be a spine of an orientable $3$-manifold $N$; in
particular its mod $2$ gleam has to be zero. By performing the
construction above, using as gleam on $P$ the $0$ gleam over all
the regions, we get the manifold $N\times [- 1,1]$.
\end{example}
\begin{rem}
All the $4$-manifolds obtained by
``thickening" the polyhedra equipped with gleams as in the proof
of Theorem \ref{teo:reconstruction} admit a handle decomposition
containing no handles of index greater than $2$. It can be shown
that also the reverse holds: any manifold admitting such a handle
decomposition can be obtained by applying Theorem
\ref{teo:reconstruction} to a suitable polyhedron equipped with
gleams (see \cite{Co}) .
\end{rem}

\begin{defi}
A polyhedron equipped with gleam $(P,gl)$ is a {\it shadow} of a
$4$-manifold $M$ if $M$ is diffeomorphic to the thickening $M_P$
of $P$ obtained through Theorem \ref{teo:reconstruction}.
\end{defi}

The above remark shows that shadows can be used to describe
combinatorially only a subset of all the smooth $4$-manifolds not
including closed ones. To obviate to this apparent weakness of the
theory, let us recall the following result due to F. Laudenbach
and V. Poenaru proved in \cite{LP}:
\begin{teo}\label{teo:LP}
Let $M$ be an oriented, smooth and compact $4$-manifold with
boundary equal to $S^3$ or to a connected sum of copies of
$S^2\times S^1$. Then there is only one (up to diffeomorphism)
closed, smooth and oriented $4$-manifold $M'$ obtained by
attaching to $M$ only $3$ and $4$-handles.
\end{teo}

Roughly speaking, the above result states that when a manifold is
``closable" then it is in a unique way (up to diffeomorphism).
This allows us to describe all the closed $4$-manifolds by means
of polyhedra with gleams: indeed given a closed manifold equipped
with an arbitrary handle decomposition, 
considering the union of all handles of index strictly less than $3$
we get a new manifold $M$ which admits a shadow and can be
described combinatorially as explained above. The initial manifold
can be then uniquely recovered from $M$ because of Theorem
\ref{teo:LP}. We can then give the following definition:
\begin{defi}
A polyhedron with gleams $(P,gl)$ is a {\it tunnel shadow} of a
closed $4$-manifold $X$ if and only if $X$ can be obtained by
attaching $3$ and $4$-handles to the $4$-manifold
$M_P$ obtained from $P$ through the reconstruction map of Theorem
\ref{teo:reconstruction}.
\end{defi}
The name ``tunnel shadow" comes from the fact that the inverse of attaching a $3$-handle is the operation of digging a tunnel along a properly embedded arc.
We will see later
ho to construct a shadow of the complement of a properly embedded arc in the thickening $M_P$ of a shadow $P$.
\begin{example}\label{cp2}
Let $P$ be equal to a $2$-sphere and let its gleam be $- 1$. 
The thickening of $(P,gl)$ is the total space of the disc bundle
over $S^2$ with Euler number is $- 1$, whose boundary is $S^3$.
This space is the complement of a point in
$\mathbb{CP}^2$, and, because of Theorem
\ref{teo:LP}, $(P,gl)$ is a tunnel shadow of $\mathbb{CP}^2$.
\end{example}

\section{Shadows and basic moves}\label{sec4}
In this section we outline the list of basic modifications
of polyhedra with gleams called {\it moves} which allow one to
produce new shadows of the same manifold from a given one.

From now on, let $(P,gl)$ be a shadow of a $4$-manifold $M$ into
which it is embedded through the construction of the
proof of Theorem \ref{teo:reconstruction}. The first three moves
we describe are called {\it shadow equivalences}; they
modify the position of $P$ in $M$ producing a
new shadow $P'$ embedded in $M$ which differs from $P$ only inside
a $4$-ball.
\begin{figure}[htbp]
  \centerline{\includegraphics[width=11.4cm]{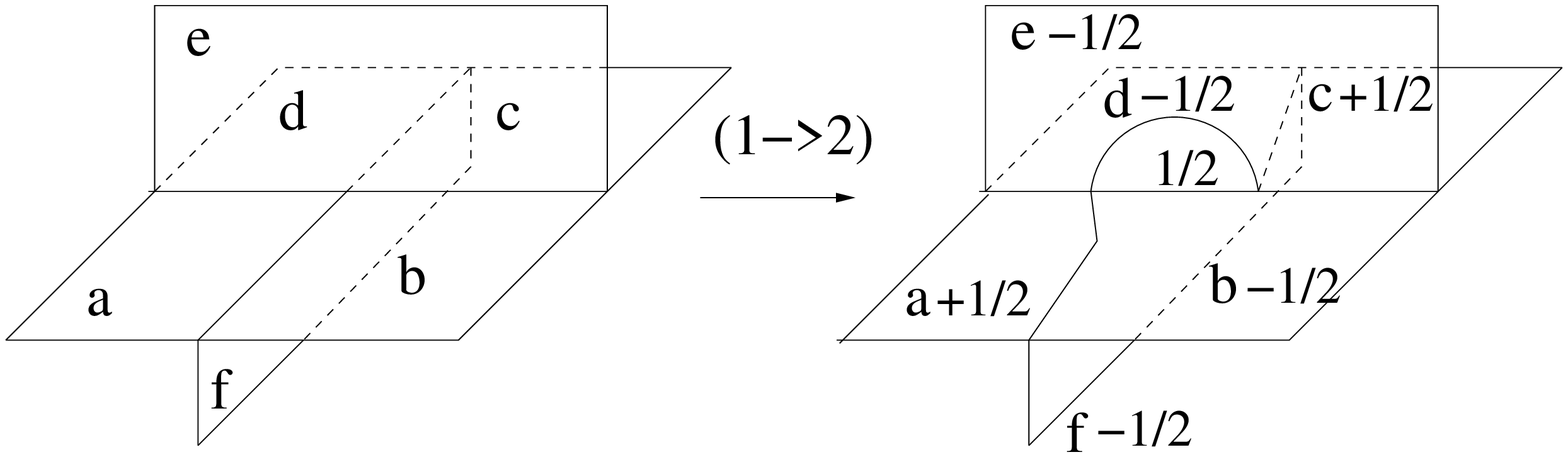}}
  \centerline{\includegraphics[width=11.4cm]{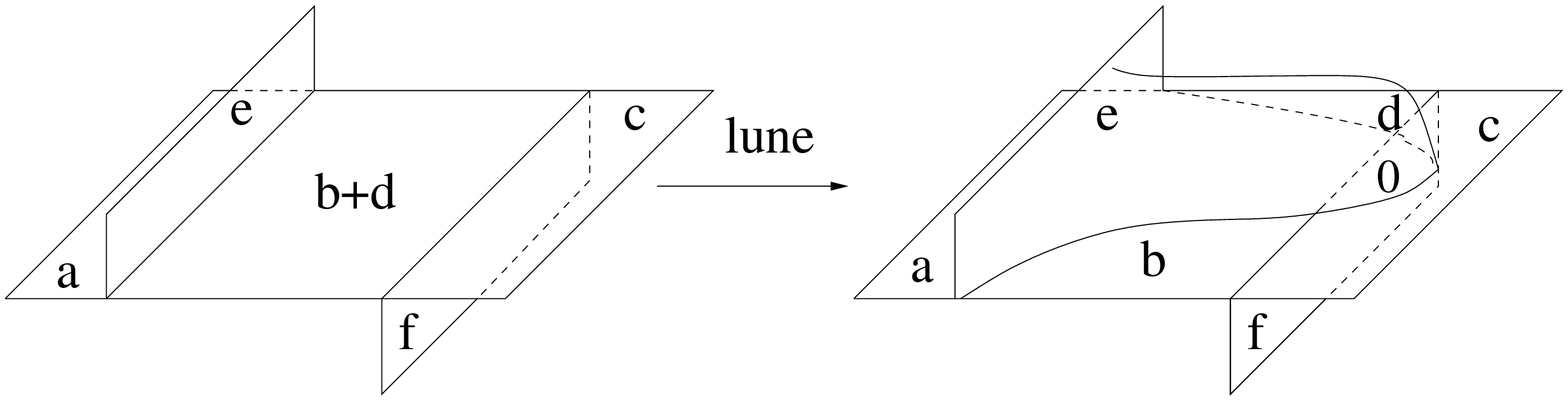}}
  \centerline{\includegraphics[width=11.4cm]{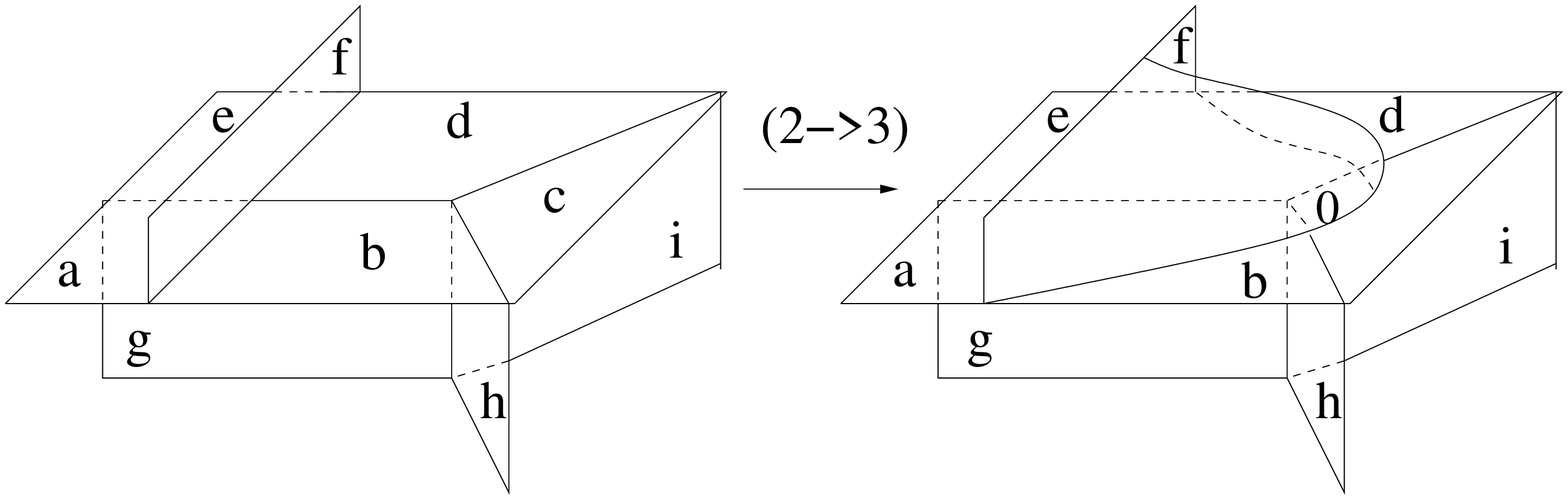}}
  \caption{The three shadow equivalences.}
  \label{fig:2-3}
\end{figure}

These moves are shown in Figure \ref{fig:2-3}; let us by the
moment forget about the letters written in the figure. To visualize these moves,
imagine that a region of $P$ is slid over
some other regions, producing a polyhedron which is different
from the initial one only in a contractile subpolyhedron (the one
pictured in the left part of the figure). These moves are called
respectively $1\to 2$, $0\to 2$ (or lune) and $2\to 3$-moves,
because of their effect on the number of vertices of the
polyhedra. Since the $1\to 2$-move is the trickiest to
understand, we concentrate first on the $0\to 2$-move (also
called finger move or lune move); this move acts in a $4$-ball
contained in $M$ and containing the part of $P$ shown in the left
part of the picture. After the move, we substitute this
part of $P$ with the new one drawn in the right part of the figure
and we keep the $4$-ball unchanged, obtaining a shadow
$P'$ of the same $4$-manifold $M$. One can imagine 
the move as a sliding over the horizontal regions of the
region which in the left part of the figure is the upper-left;
this move can be performed in a
$3$-dimensional slice of the $4$-ball in $M$; hence we say that the $0\to 2$-move
has a $3$-dimensional nature. The same holds for the $2\to
3$-move, also called Matveev-Piergallini move: although
apparently more complicated than the $0\to 2$-move, it can
be seen that it is a version of the $0\to
2$-move in which the sliding region passes over a vertex. The
same comments as above apply to this case.

The case of the $1\to 2$-move is different; this
move applies to a neighborhood of a vertex and lets a region (the
vertical lower one in the left part of the figure) slide over the
vertex producing a new one. It is a good exercise to try to
visualize this sliding: it cannot be performed in $\mr^3$; if we
want the sliding region not to touch in its interior the
other regions near the vertex, then it is necessary to use the
fourth dimension. Hence we say that this move has $4$-dimensional
nature; it has the effect of
modifying $P$ but not $M$. Note that, even if  the move ``needs
the fourth dimension", its result is still a locally flat
polyhedron in $M$.

Of course, the above comments apply also to the inverses of the
moves, that is to the moves which in the figure go from right to left.
Let us now clarify the meaning of the letters written on the
regions of the polyhedra in the figure: they
represent the gleams of the regions; note that only the
$1\to 2$-move changes the gleams of the regions near the vertex
over which it is applied. Each
move represents a modification of the embedded polyhedron, and
hence produces a new polyhedron whose gleams (induced by the
embedding in the ambient manifold as explained in Proposition
\ref{embpoly}) can differ from those of the initial
polyhedron; moreover these gleams can change only on the regions
which appear in the pictures, since
otherwise the position of the regions in $M$ is left unchanged by
the move.

All these comments can be resumed in the following proposition:
\begin{prop}
Let $P$ be a shadow of a $4$-manifold $M$ and let $P'$ be obtained
from $P$ through the application of any sequence of $1\to 2$,
$0\to 2$ and $2\to 3$ embedded moves (and possibly their
inverses), then $P'$ is a shadow of $M$. Moreover, the gleam of
$P'$ can be reconstructed from the gleam of $P$ by following step
by step the changes on the gleams induced by the moves of the
sequence.
\end{prop}

\begin{figure}[htbp]
  \centerline{\includegraphics[width=11.4cm]{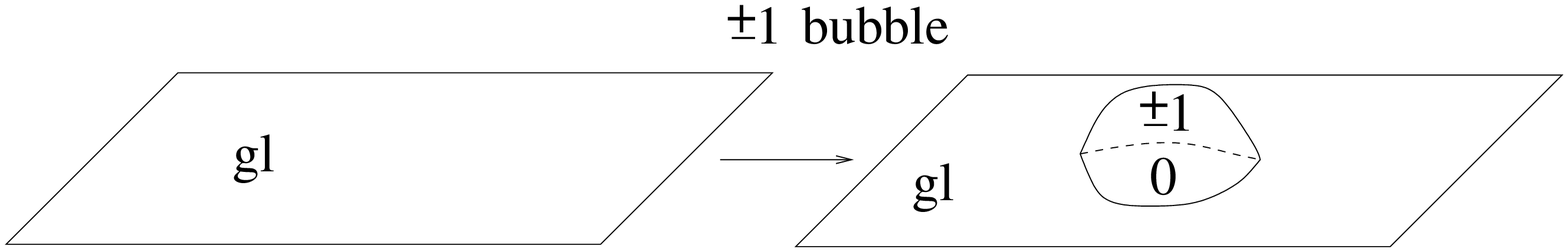}}
  \caption{The $\pm 1$-bubble moves: the gleam of the region over which the move is applied is left unchanged.
  A new $0$-gleam disc is created and a $\pm 1$-gleam disc is attached to the boundary of the latter disc.
  The $0$-bubble move, creates a $0$-gleam disc in place of the $\pm 1$-gleam one.}
  \label{fig:bubble}
\end{figure}

The above moves change the polyhedron but do not change
the ambient manifold; we now describe two moves which
change also the manifold and which are useful in a number of
applications we will outline in the next sections: 
the $0$-bubble move and the $\pm 1$-bubble move. On the
level of naked polyhedra they are identical and shown in Figure
\ref{fig:bubble}; they differ on the level of the gleams of the
polyhedron after the move. The $0$-bubble move applied to a polyhedron $P$ 
produces a polyhedron $P'$ by attaching a $0$-gleam disc to a
$0$-gleam disc contained in a region of $P$ and
leaving unchanged the gleams of the other regions of the
polyhedron; $P'$ is not standard, but
the construction of Theorem \ref{teo:reconstruction} can be generalized to simple polyhedra equipped with gleams. Let $M_P$ and $M_{P'}$ be the $4$-manifolds determined
respectively by $P$ and $P'$ as explained in Theorem
\ref{teo:reconstruction}; it can be checked that $M_{P'}$ is
diffeomorphic to the complement in $M_P$ of a properly embedded
arc, that is, it is obtained
by tunneling $M_P$. Conversely, $M_P$ can be recovered from
$M_{P'}$ by attaching a $3$-handle whose cocore corresponds
to that arc.

Let us now describe the effects of the $\pm 1$-bubble moves. These
moves are identical except the sign chosen; we
concentrate on the $- 1$-bubble move. Its application to a polyhedron $P$ 
produces a polyhedron with gleams $P'$ by attaching a $- 1$-gleam disc
over a $0$-gleam disc in a region of $P$ and leaving the other
gleams unchanged. It can be checked that if $M_P$ is the manifold
determined by $P$ then the one determined by $P'$ is
$M_{P'}=M_P \#\mathbb{CP}^2$, and, in particular, $\partial
M_{P}=\partial M_{P'}$; analogously the application of a
$1$-bubble move produces a shadow of
$M_P\#\overline{\mathbb{CP}}^2$.

\begin{figure}[htbp]

  \centerline{\includegraphics[width=11.4cm]{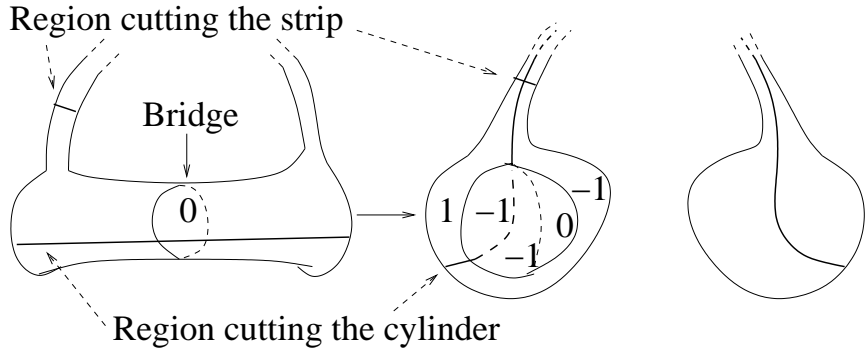}}

  \caption{In the left part of the figure we show the pattern to
  which the trading-move applies. It is constituted of a cylinder cut by a
  zero gleam disc and of a (possibly twisted) strip connecting the two ends of the cylinder. Some strands (corresponding
  to attaching curves of regions of the shadow) can pass over
  the cylinder or cut the strip. In the right part we show the effect
  of the trading-move: the cylinder disappears, its two ends are
  capped off with two zero gleam discs and, in one of them, a little $0$-gleam two
  disc is attached. The strands which before the move run
  across the cylinder, are now completed along the strip and linked
  with the little zero gleam disc as shown in the drawing.}
  \label{fig:trading}
\end{figure}
The last move we introduce is called the ``trading''-move and 
is visible in Figure \ref{fig:trading}. Unlike the above moves, 
this move does not start from a contractible
sub-polyhedron of the initial polyhedron; hence we say that it is
``non-local'': one can apply it only when a pattern in the
polyhedron which coincides with the left part of the figure is
found. The move dramatically changes the topology of
the polyhedron: in particular the fundamental group changes. For
the sake of simplicity, in the figure we draw a version of the
trading move in the case of a non standard, simple polyhedron;
up to the application of a suitable
number of $0\to 2$-moves, the polyhedron can be transformed into a
standard one and the move can be described also in the standard
setting but the drawings are more complicated and we omit them here.

\begin{figure}[htbp]
  \centerline{\includegraphics[width=7.4cm]{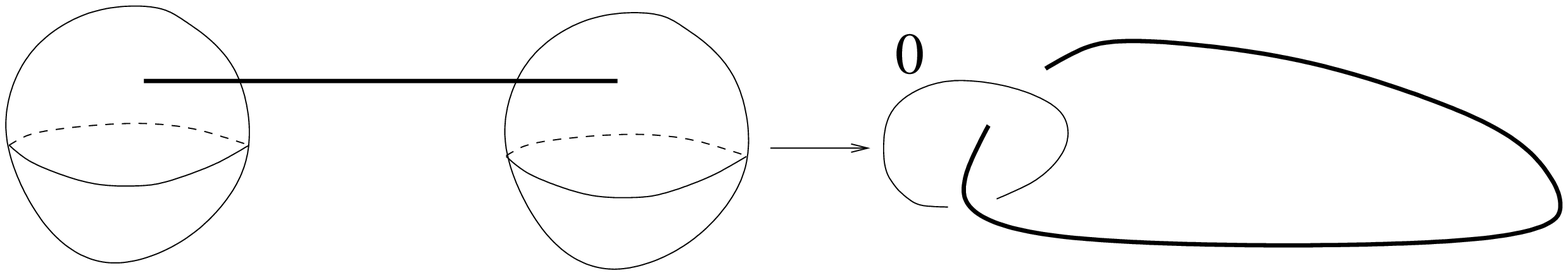}}
  \caption{In the left part of the figure, we show an example of $1$-handle in a Kirby
  calculus presentation of a $4$-manifold over which an attaching curve of a $2$-handle passes.
  In the right part we show
  the effect of trading the $1$-handle for a $2$-handle.}
  \label{fig:kirbytrade}
\end{figure}

\begin{prop}[\cite{Tu}]\label{lem:trademove}
Let $P$ be an integer shadowed polyhedron and let $P'$ be the
polyhedron obtained by applying a trading-move (or its inverse) to
$P$. Then the manifold $M_{P'}$ is obtained from $M_{P}$ by doing
surgery along a closed curve $C$ embedded in the interior of
$M_{P}$ (respectively, on an embedded sphere whose self
intersection number is zero), and then, in particular, $\partial
M_{P}$ and $\partial M_{P'}$ are diffeomorphic.
\end{prop}
\begin{prf}{1}{ We limit ourselves to rough motivations for the above proposition; see 
\cite{Tu} for details. We first
analyze the portion of manifold described in Figure
\ref{fig:trading} before the move. We call {\it bridge} the part
of the polyhedron $P$ formed by the curved cylinder together with
the central disc with zero gleam in the left part of Figure
\ref{fig:trading}. It can be checked that there exists a handle
decomposition of $M_P$ such that the bridge is contained in a
$1$-handle and the two ends of the bridge are contained in the two
ends of the $1$-handle. The other regions which in the figure run
along the cylinder correspond to $2$-handles whose attaching
curves pass through the $1$-handle. So, from the point of view of
the Kirby calculus, the left part of Figure \ref{fig:trading}
corresponds to the left part of Figure \ref{fig:kirbytrade}, where
the strand running from one attaching sphere of the $1$-handle to
the other one corresponds to the attaching curve which in Figure
\ref{fig:trading} runs across the bridge.

Now, let $C$ be a closed curve inside $M_P$ which passes once over the
bridge in Figure \ref{fig:trading} and then gets
closed the strip connecting the two ends of the cylinder;
surgering along $C$ produces a $4$-manifold which, in
place of the one-handle determined by the bridge, contains a 
two-handle whose attaching curve is an unknot in $\partial
M_P$ and such that the attaching curve of every $2$-handle which
passed over the bridge, passes inside the disc bounded by the
unknot (recall that surgering along a curve in a $4$-manifold
consists of substituting the regular neighborhood of the curve
with a copy of $S^2\times D^2$).

It can be checked, but we will not do it
here, that this move corresponds to the move which in Kirby
calculus, given a presentation of a $3$-manifold, allows one to
exchange a one handle with a zero framed two handle represented by
an unknot and linked with all the two handles passing over
the one handle, obtaining a presentation of the same $3$-manifold
(see the right part of Figure \ref{fig:kirbytrade}).

}\end{prf}

We resume the facts exposed in this section in the following
theorem. For a complete account we refer to \cite{Tu}.
\begin{teo}\label{teo:moveeffect}
Let $P$ and $P'$ be two integer shadowed polyhedra such that $P'$
is obtained from $P$ by applying a move of the following types:
\begin{enumerate}
\item $2\to 3$-move or its inverse;
\item finger-move or its
inverse;
\item $1\to 2$-move or its inverse.
\end{enumerate}
Then $M_P$ and $M_{P'}$ are diffeomorphic by an orientation
preserving diffeomorphism.

If $P'$ is obtained from $P$ by a $\pm 1$-bubble move then
$M_{P'}$ is diffeomorphic to the boundary connected sum of $M_P$
and a punctured $\mathbb{CP}^2$ or $\overline{\mathbb{CP}^2}$
respectively, hence $\partial M_{P'}$ is diffeomorphic (by an
orientation preserving diffeomorphism) to $\partial M_P$. If $P'$
is obtained from $P$ by a $0$-bubble move then $M_{P'}$ is
diffeomorphic to the boundary connected sum of $M_P$ and
$S^2\times D^2$, and then, in particular $\partial M_{P'}=\partial
M_P \# S^2\times S^1$. Finally, if $P'$ is obtained from $P$ by
the application of a trading move, then $M_{P'}$ is obtained from
$M_P$ by doing a surgery along the curve which passes once over
the bridge and once over the cylinder of Figure \ref{fig:trading};
in particular $\partial M_{P'}=\partial M_P$.

\end{teo}

\section{Some applications and questions}\label{sec5}
In this section we briefly recall some results and open questions. 
The first one surges once one proves
Theorem \ref{teo:reconstruction} and is the following:
\begin{question}
Given two shadows $P$ and $P'$ of the same manifold, can they be
connected by a sequence of suitable moves, for instance shadow
equivalences (see the beginning of Section \ref{sec4})?
\end{question}

A first answer to the above question is due to Turaev \cite{Tu}:
\begin{teo}\label{teo:calculus}
Let $P$ and $P'$ be two shadows of the same $4$-manifold $M$, then
they can be connected by a sequence of shadow equivalences,
$0$-bubble moves and their inverses. Moreover the sequence always
contains a number of $0$-bubbles equal to the number of inverse
$0$-bubble moves.
\end{teo}
In the above result there is a disturbing aspect: the sequence
of moves connecting $P$ and $P'$ can contain $0$-bubble moves
and can produce at a certain step a polyhedron which is
not a shadow of $M$, hence in general the sequence does not pass only through
shadows of $M$. Moreover the converse of the above theorem is not
yet proved:
\begin{question}
Let $P$ and $P'$ be obtained from each other by a sequence of
shadow equivalences and $0$-bubble moves containing a number of
$0$-bubble moves equal to the number of their inverses. Is it true
that the manifold $M$ and $M'$ determined by $P$ and $P'$
respectively are diffeomorphic?
\end{question}

The author in \cite{Co} proved that the answer to the above
question is ``yes" when $P$ is (and hence $P'$) is simply
connected; the positive answer in general is still a conjecture.

The following question can be thought of as an
embedded version of the Andrews-Curtis conjecture:
\begin{question}
Let $P$ and $P'$ be shadows of the same manifold; can one connect
them by a sequence of shadow equivalences? (No bubble moves!)
\end{question}
By Matveev and Piergallini's results  (\cite{Ma},
\cite{Pi}), the analogous question in dimension $3$ has a
positive answer. The $4$-dimensional case remains open;
it has to be pointed out here that, although apparently the above
question is strictly related to the Andrews-Curtis conjecture, by
the moment no formal relation between the two problems has been
proved.

Theorem \ref{teo:calculus} can be re-read in terms of tunnel
shadows of closed $4$-manifolds, and it then states that any two
shadows of a closed $4$-manifold are connected by a sequence of
shadow equivalences and $0$-bubble moves. In this case, it
represents a genuine calculus in the standard sense.

\section{Shadows of $3$-manifolds}\label{sec6}
In this section we will briefly recall the notion of a shadow of a $3$-manifold and give an overview of some recent results on the
subject.

\begin{defi}
Let $(P,gl)$ be a shadow of a $4$-manifold $M$; then $(P,gl)$ is a shadow of the three-manifold $N=\partial M$.
\end{defi}

It can be proved that any closed and oriented $3$-manifold admits
a shadow, and even a simply connected one; moreover, the
application of any shadow equivalence or $\pm 1 $-bubble move to a
shadow of a $3$-manifold produces a shadow of the same
manifold. In general, the following was proved in \cite{CT}:
\begin{teo}[Calculus for shadows of $3$-manifolds]\label{teo:3calculus}
Any two simply connected shadows of the same closed and oriented
$3$-manifold can be connected by a suitable sequence of shadow
equivalences, $\pm 1 $-bubble moves and their inverses.
\end{teo}

The above result is an analogue of the Kirby calculus in
the world of shadows of $3$-manifolds, where the blow-up move is
substituted by the $\pm 1$-bubble moves; its pretty feature is
that it is based only on local moves.

A new notion of complexity of $3$-manifolds based on shadows is
defined in \cite{CT2}:
\begin{defi}
Let $N$ be a closed and oriented $3$-manifold. The {\it shadow
complexity} of $N$ is the minimal number of vertices of a
(possibly non standard) shadow of $N$.
\end{defi}
\begin{example}
As shown in Example \ref{cp2}, the polyhedron $S^2$ equipped with
gleam $- 1$ is a shadow of $S^3$ and hence the shadow complexity of
$S^3$ is zero. \end{example}

Strong relations have been found between shadow complexity and
hyperbolic geometry of $3$-manifolds. The following Theorem,
to appear in \cite{CT2}, summarizes some of them:
\begin{teo}
There exist universal positive constants $h_1$ and $h_2$ such that
for any hyperbolic, complete $3$-manifold $N$ with finite volume
the following inequalities hold: $h_1 Vol(N)\leq sc(N)\leq h_2
(Vol(N))^2$, where $sc(N)$ is the Shadow Complexity of $N$ and
$Vol (N)$ is its hyperbolic volume.
\end{teo}
Moreover, in the same paper we prove that each $3$-dimensional graph
manifold admits a shadow containing no
vertices.

These results suggest that a study of hyperbolic manifolds by
means of shadows can be attempted; to that purpose, one of the
first questions to be answered is the following:
\begin{question}
Given a polyhedron with gleams $(P,gl)$, what are the necessary
and sufficient conditions on $(P,gl)$ ensuring that the
$3$-manifold determined by the shadow is irreducible, or atoroidal
or hyperbolic?
\end{question}
A clear and definitive combinatorial condition answering the above
question is still missing. Partial answers have been obtained in
\cite{CFMP}: for instance, it has been proved that if a
$3$-manifold admits a standard shadow whose gleams are all greater
than $7$ in absolute value, then the manifold admits a metric with
negative sectional curvature, and, consequently, is aspherical and atoroidal.

\section{A glance on branched shadows}\label{sec7}
Besides the open problems exposed in the preceding sections,
which, in part, are internal to the theory of shadows, a study of
additional structures and invariants of $4$-manifolds can be
attempted by means of shadows. In this section we 
summarize the results proved by the author in \cite{Co} (and which
will be exposed in a self contained way in \cite{Co3} and
\cite{Co4}) in this direction and based on the notion of branched
shadow.

Given a simple polyhedron $P$ we define the notion of {\it
branching} on it as follows:
\begin{defi}[Branching]\label{branchingcondition}
A branching $b$ on $P$ is a choice of an orientation for each
region of $P$ such that for each edge of $Sing(P)$, the
orientations induced on the edge by the regions containing it do
not coincide.
\end{defi}

\begin{figure} [h!]
   \centerline{\includegraphics[width=11.4cm]{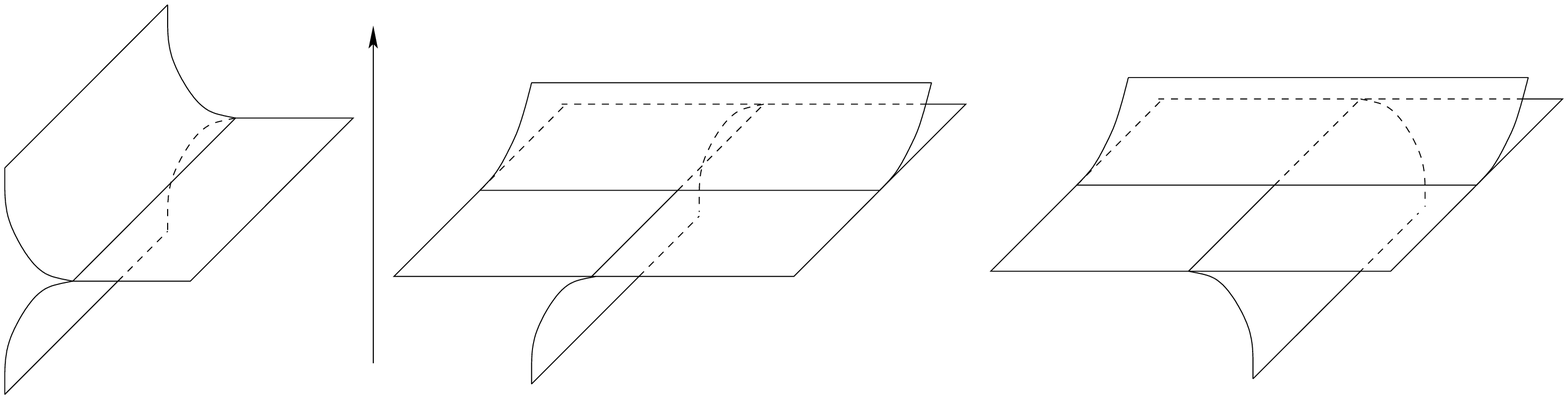}}
  \caption{How a branching allows a smoothing of the polyhedron.}\label{branching}
\end{figure}

Not all the simple polyhedra admit a branching and, on the
contrary, there are some which admit many different ones. So, we
will say that a polyhedron is {\it branchable} if it admits a
branching and we will call {\it branched polyhedron} a pair
$(P,b)$ where $b$ is a branching on the polyhedron $P$.
\begin{defi}
Let $(P,gl)$ be a shadow of an oriented $4$-manifold $M$; $P$ is
said to be {\it branchable} if the underlying polyhedron is. We
call {\it branched shadow} of $M$ the triple $(P,gl,b)$ where
$(P,gl)$ is a shadow and $b$ is a branching on the underlying
polyhedron. When this does not cause any confusion, we do not
specify the branching $b$ nor the gleam $gl$ and we will simply
write $P$.
\end{defi}
A branching on a simple polyhedron allows us to smoothen its
singularities and equip it with a smooth structure as shown in
Figure \ref{branching}.

Using the orientation and the smooth structure of the shadow, in
\cite{Co3} (see also \cite{Co}), we will show how to equip the regular neighborhood of a
shadow in a $4$-manifold with a pair mutually transverse
distributions of oriented $2$-planes (roughly speaking, the planes
which in Figure \ref{branching} are horizontal and their
orthogonal complements). Then, in the same paper using this pair of distributions of
oriented $2$-planes, we show how to fix a homotopy class of almost
complex structures on the regular neighborhood of the shadow, and
hence a $Spin^c$-structure on this manifold. The
result to appear in \cite{Co3} can be summarized as follows:
\begin{teo}\label{teo:spinc}
Let $M$ be a $4$-manifold admitting a shadow and let $s$ be a
$Spin^c$-structure on $M$. There exists a branched shadow $P$ of
$M$ representing the $Spin^c$-structure $s$. Stated
differently, given a homotopy class of almost complex structures
$[J]$ on $M$, there exists a branched shadow of $M$ representing
the homotopy class $[J]$ of almost complex structures on $M$.
\end{teo}

Following the line of the last part of Theorem
\ref{teo:spinc}, it can be studied how branched shadows can be
used to find genuine complex structures on their thickening. In
\cite{Co3}, the problem is attacked through techniques that are similar
to those already used for real surfaces in complex manifolds by
Bishop in \cite{Bi} and Harlamov-Eliashberg in \cite{HE}, based on
counting the number of complex elliptic and hyperbolic points. The
result we obtain can be summarized as follows:
\begin{teo}[\cite{Co3}]\label{teo:shadowcondition}
Let $P$ be a branched shadow reconstructing the manifold $M$
equipped with the almost complex structure $J$. There exists a
genuine complex structure $J'$ homotopic to $J$ such that $P$
contains no $J'$-complex points of negative type, that is, $P$
contains no points where its tangent plane is $J'$-complex but
oriented differently by the branching on $P$ and its complex structure induced by $J'$.
\end{teo}

As a natural development of the above result, using
Gompf-Eliashberg's constructions, combinatorial sufficient
conditions on a branched shadow ensuring that its thickening
admits a Stein domain structure will be studied in \cite{Co4} (see also \cite{Co}, Chapter V ). The conditions we will find are based on the definition (which we
will not recall here) for each region of a branched polyhedron of
an integer coefficient called the Euler number which is easily
calculable from the combinatorial structure of the polyhedron. The
main result, which will be exposed in a self-contained way in \cite{Co4}, is
summarized in a simplified version by the following:
\begin{teo}[\cite{Co4}]
Let $P$ be a branched shadow of a four manifold $M$. There exist constants $k_i$ such that if, for
each region $R_i$ of $P$ it holds $Eul_i+gl_i\leq k_i$ where $gl_i$
is the gleam of $R_i$ and $Eul_i$ is its Euler number, then $M$
admits a Stein structure whose underlying complex structure
belongs to the homotopy class of almost complex structures
determined by $P$.
\end{teo}

The conditions found in the above theorem have been exploited in \cite{Co}
to explicitly exhibit minimal genus representatives of some
particular homology classes of surfaces in $4$-manifolds:

\begin{cor}\label{cor:oertel4}
Let $P$ be a branched shadow satisfying the hypotheses of Theorem
\ref{teo:shadowcondition} and such that $gl_i\geq 0$ for every
$i$, and let $[S]\in H_2(P;\mathbb{Z})$ be a homology class of the form $\sum_i n_i R_i$ where $R_i$ are the regions of $P$ and $n_i$ are non negative integers. Then
there is an explicit construction of an oriented embedded
representative $S'\subset M$ of $[S]$ having the lowest possible genus in
its homology class.
\end{cor}

Moreover, using Theorem \ref{teo:shadowcondition} on branched
shadows equipped with zero gleams on every region, the following
result, which was conjectured by Benedetti and Petronio in
\cite{BP} is proved:
\begin{teo}
Let $P$ be a branched spine of an oriented $3$-manifold $N$ such
that the Euler number of each region of $P$ is non positive, then
the homotopy class of oriented $2$-planes encoded by the branching
of the spine on its $3$-dimensional thickening contains a
representative distribution which is a tight contact structure.
\end{teo}

The above results are likely to be non optimal and
further study is required to fully understand the potential of the
theory of branched shadows.

\noindent

\end{document}